\documentclass[reqno]{amsart}
\usepackage{amssymb}
\usepackage{amsaddr}
\usepackage{amscd}

\newtheorem{theorem}{Theorem}[section]

\newtheorem{lemma}[theorem]{Lemma}
\newtheorem{proposition}[theorem]{Proposition}

\newtheorem{remark}{Remark}[section]

\newtheorem{definition}{Definition}[section]
\numberwithin{equation}{section}

\begin{document}
\title[Nonlinear\ Schr\"{o}dinger equation]{Small data blow-up of $L^2$-solution
 for the nonlinear Schr\"{o}dinger equation without gauge invariance}
\author[M. Ikeda and Y. Wakasugi]{Masahiro Ikeda$^1$ and Yuta Wakasugi$^2$}
\address{{\rm Department of Mathematics
Graduate School of Science, Osaka University,\\
Toyonaka Osaka, 560-0043, Japan\\[5pt]
{\footnotesize E-mail addresses : $^1$m-ikeda@cr.math.sci.osaka-u.ac.jp\\
\hspace{10em}$^2$y-wakasugi@cr.math.sci.osaka-u.ac.jp
}}}
\subjclass[2000]{35Q55}


\begin{abstract}
We study the initial value problem for the nonlinear Schr\"odinger equation
$$
	i\partial _{t}u+\Delta u=\lambda\vert u\vert ^{p},
	\quad\left( t,x\right) \in \left[ 0,T\right) \times \mathbf{R}^{n},
$$
where $1<p\leq 1+2/n$ and $\lambda\in\mathbf{C}\setminus\{0\}$.
We will prove that the blow-up of the $L^2$-norm of solutions with suitable initial data.
We impose a condition related to the sign of the data but put no restriction on their size.
\end{abstract}
\maketitle

\section{\label{S1} Introduction}

We study existence of a blow-up solution for the nonlinear Schr\"{o}dinger
equation (NLS) with a critical or subcritical non-gauge invariant power type
nonlinearity 
\begin{equation}
	i\partial _{t}u+\Delta u=\lambda \left\vert u\right\vert ^{p},\quad
	\left(t,x\right) \in \left[ 0,T\right) \mathbf{\times R}^{n},\label{eq11}
\end{equation}
with the initial condition
\begin{equation}
	u\left( 0,x\right) =f\left( x\right) ,\text{ \ }x\in \mathbf{R}^{n},
\label{eq12}
\end{equation}
where
$T>0,$ $1<p\leq 1+2/n,$ $u=u\left( t,x\right)$
is a complex-valued unknown function,
$\lambda=\lambda _{1}+i\lambda _{2}\in \mathbf{C}\setminus\left\{ 0\right\}$,
$\lambda _{j}\in \mathbf{R}$ $\left( j=1,2\right) ,$
$f=f\left( x\right)=f_{1}\left( x\right) +if_{2}\left( x\right)$
and
$f_{j}=f_{j}\left( x\right)\in L^1_{loc}(\mathbf{R}^n)\, \left( j=1,2\right)$
are real-valued functions.

It is well known that local well-posedness holds for (\ref{eq11}) in
several Sobolev spaces $H^{s}$ $\left( s\geq 0\right) $
(see e.g. \cite{Caze03, Tsu87} and the references therein).
However, there is no result of global existence of the solution of (\ref{eq11})-(\ref{eq12}).
In this paper, we will prove that
if the initial data $f$ in $L^2$ satisfies a certain condition related to its sign,
then the $L^2$-norm of the solution $u$ of (\ref{eq11})-(\ref{eq12}) blows up in finite time,
even if the data is sufficiently small (see section 2).
We note that when $p\geq p_s$, where $p_s$ is the well-known Strauss exponent,
which is greater than $1+2/n$, global existence results are known (see \cite{Caze03}).
However, it is still open what happens in the case of $1+2/n<p\leq p_{s}$.

Our result implies that the nonlinear effect of
$\lambda \left\vert u\right\vert ^{p}$
is quite different from that of
$\lambda _{0}\left\vert u\right\vert ^{p-1}u$
$\left( \lambda_{0}\in \mathbf{R}\right)$,
since the $L^{2}$-norm of solutions for 
\begin{equation}
i\partial _{t}u+\Delta u=\lambda _{0}\left\vert u\right\vert ^{p-1}u\label{1.10}
\end{equation}%
conserves for any $t\in \mathbf{R}.$
Tsutsumi \cite{Tsu87} proved global existence of $L^{2}$-solution
of (\ref{1.10}) when $1<p<1+4/n$.
It is also well known that for (\ref{1.10}), the exponent $p=1+2/n$ is the threshold between
the short range scattering and the long range one
(see \cite{Ba84, TsuYaji, Oza91, Gi-Oza93, HaNa98, HaKaNa99}).
We also mention that when $p\geq 1+4/n$, blow-up of $H^1$-solution of (\ref{1.10}) is proved by
Glassey \cite{Gla77} (see also \cite{OgaTsu91}).
However, their results require that the data are large as contrast with our result.

Back to our problem (\ref{eq11}),
in the critical case $\left( n,p\right) =\left( 2,2\right)$,
Shimomura \cite{Shi05} and Shimomura-Tsutsumi \cite{Shi-Tsu06}
studied the asymptotic behavior of solutions of (\ref{eq11}).
Especially, Shimomura-Tsutsumi \cite{Shi-Tsu06}
proved nonexistence of the wave operator for (\ref{eq11}).
On the other hand, Hayashi-Naumkin \cite{HaNa06}
considered the final state problem for NLS with the quadratic nonlinearity
$\mu u^2+\nu \bar{u}^2+\lambda |u|^2$,
which includes the term
$\lambda |u|^2$,
in two space dimension.
They proved existence of the global solution which behaves unlike the free one in $L^2$.
We note that their result requires that $\mu, \nu \neq 0$ and
is not applicable to (\ref{eq11}).

From these results,
some people might think that the critical or subcritical non-gauge invariant
nonlinearity $\lambda \left\vert u\right\vert ^{p}$ may act as a long range
effect such as $\lambda _{0}\left\vert u\right\vert ^{p-1}u$.
However, our result gives a negative conclusion to such an expectation.

\section{\label{S2} Main Result}

We first recall the well-known fact about local existence of the solution in $L^2$ for the integral equation
\begin{equation}
u\left( t\right) =U\left( t\right) f-i\lambda \int_{0}^{t}U\left( t-s\right)
\left\vert u\right\vert ^{p}ds  \label{eq21}
\end{equation}
associated with (\ref{eq11})-(\ref{eq12}),
where $U(t)=\exp(it\Delta)$ is the evolution group of the free Schr\"odinger operator.

\begin{proposition}[Tsutsumi \cite{Tsu87}]\label{Prop21}
Let
$1<p<1+4/n, \lambda \in \mathbf{C}\setminus\left\{ 0\right\}$
and
$f\in L^{2}$.
Then there exist a positive time
$T=T\left(\left\Vert f\right\Vert _{L^{2}}\right) >0$
and a unique solution
$u\in C\left( \left[ 0,T\right) ;L^{2}\right)\cap
L_{t}^{r}\left( 0,T ;L_{x}^{\rho }\right)$
of the integral equation {\rm (\ref{eq21})},
where $r, \rho$ are defined by
$\rho =p+1$ and $2/r=n/2-n/\rho$.
\end{proposition}

We call the solution $u$ in the above proposition ``$L^2$-solution".
Let $T_m$ be the maximal existence time of local $L^{2}$-solution, that is,
\begin{eqnarray*}
	T_{m} &\equiv &\sup \left\{ T\in \left( 0,\infty \right] ;
		\text{ there exists the unique solution }u\text{ to (\ref{eq21})}\right. \\
	&&\left. \text{such that }u\in
	C\left( \left[ 0,T\right) ;L^{2}\right) \cap
	L_{t}^{r}\left( 0,T ;L_{x}^{\rho }\right) \right\},
\end{eqnarray*}%
where $r, \rho$ are as in the above proposition.
To state our result, we put the following assumption on the data:
\begin{equation}
\label{eq22}
	``f_1\in L^1(\mathbf{R}^n), \ \lambda_2\int_{\mathbf{R}^n}f_1(x)dx>0"\ \text{or}\ 
	``f_2\in L^1(\mathbf{R}^n), \ \lambda_1\int_{\mathbf{R}^n}f_2(x)dx<0".
\end{equation}
Our main result is the following:

\begin{theorem}\label{Thm22}
Let $1<p\leq 1+2/n$,
$\lambda\in\mathbf{C}\setminus\left\{ 0\right\}$
and
$f\in L^{2}$.
If the initial data $f$ satisfies {\rm (\ref{eq22})},
then $T_{m}$ must be finite. Moreover, we have
\begin{equation}
	\lim_{t\rightarrow T_{m}-0}\left\Vert u\left( t\right) \right\Vert _{L^{2}}=+\infty.  \label{eq23}
\end{equation}
\end{theorem}
We note that we put no restriction on the size of the data.
In order to prove Theorem \ref{Thm22}, in the next section,
we introduce a weak solution of (\ref{eq11})-(\ref{eq12})
and the result of nonexistence of a global weak solution.

\section{Reduction of the problem}
To prove Theorem \ref{Thm22}, we define a weak solution of (\ref{eq11})-(\ref{eq12}).
\begin{definition}
\label{Def31}
Let $T>0$.
We mean $u$ is a weak solution of {\rm NLS} {\rm (\ref{eq11})-(\ref{eq12})} on $[0,T)$
if $u$ belongs to
$L_{loc}^p([0,T) \times \mathbf{R}^{n})$
and satisfies
\begin{eqnarray}
	\lefteqn{\int_{[ 0,T)\times\mathbf{R}^{n}}u( -i\partial _{t}\psi+\Delta\psi)dxdt}\notag\\
	&&=i\int_{\mathbf{R}^{n}}f\left( x\right) \psi \left( 0,x\right) dx+\lambda
	\int_{\left[ 0,T\right) \times \mathbf{R}^{n}}\left\vert u\right\vert^{p}\psi dxdt \label{eq31} 
\end{eqnarray}%
for any $\psi \in C_{0}^{2}\left( \left[ 0,T\right) \times \mathbf{R%
}^{n}\right) .$ Moreover, if $T>0$ can be chosen as any positive number, $u$
is called a global weak solution for {\rm (\ref{eq11})-(\ref{eq12})}.
\end{definition}

We note that an $L^{2}$-solution as in Proposition \ref{Prop21} is always a weak solution in the
sense of Definition \ref{Def31}:

\begin{proposition}
\label{Prop31} Let $T>0.$ If $u$ is an $L^{2}$-solution for the
equation {\rm (\ref{eq21})} on $\left[ 0,T\right) ,$ then $u$ is also a weak
solution on $[0,T)$ in the sense of {\rm Definition \ref{Def31}}.
\end{proposition}
We will give a proof of this proposition in Appendix.

Next, we mention nonexistence of a nontrivial global weak solution for (\ref{eq11})-(\ref{eq12})
with the condition (\ref{eq22}).

\begin{proposition}
\label{Prop32}
Let $1<p\leq 1+2/n$,
$\lambda\in\mathbf{C\backslash }\left\{ 0\right\} $
and
let $f$ satisfy {\rm (\ref{eq22})}.
If there exists a global weak solution $u$ of {\rm (\ref{eq11})-(\ref{eq12})},
then $u=0$.
\end{proposition}

Combining Proposition \ref{Prop31} and \ref{Prop32},
we obtain Theorem \ref{Thm22}.
Indeed, let $f\in L^2$ satisfy (\ref{eq22})
and $u$ be the $L^2$-solution of (\ref{eq21}).
Suppose that $T_m=\infty$.
By Proposition \ref{Prop31}, $u$ is also a global weak solution of (\ref{eq11})-(\ref{eq12})
in the sense of Definition \ref{Def31}.
Thus, we can apply Proposition \ref{Prop32}
and have $u=0$.
However, by noting $u\in C([0,\infty);L^2(\mathbf{R}^n))$,
it contradicts $f\neq 0$.
Therefore, we have $T_m<\infty$.

Next, we prove (\ref{eq23}).
First we suppose
\begin{equation*}
	\liminf_{t\rightarrow T_{m}-0}\Vert u( t) \Vert _{L^{2}}<\infty .
\end{equation*}
Then there exist a sequence
$\{ t_{k}\} _{k\in\mathbf{N}}\subset [0,T_{m})$
and a positive constant $M>0$ such that
\begin{gather}
	\lim_{k\to\infty }t_{k}=T_{m}  \label{eq32} \\
	\sup_{k\in \mathbf{N}}\Vert u( t_{k}) \Vert _{L^{2}}\leq M.  \label{eq33}
\end{gather}
By (\ref{eq33}) and Proposition \ref{Prop32}, there exists a positive constant
$T\left( M\right) $ such
that we can construct a solution%
\begin{equation*}
	u\in C( [ t_k,t_k+T( M) ) ;L^{2}) \cap L_{t}^{r}( [ t_k,t_k+T(M) ) ;L_{x}^{\rho })
\end{equation*}%
of (\ref{eq21}) for all $k\in\mathbf{N}$.
However, by (\ref{eq32}), when $k$ is sufficiently large, the inequality
$t_k+T\left( M\right) >T_{m}$
holds and it contradicts the definition of $T_{m}$.
Therefore, we obtain
\begin{equation*}
\liminf_{t\rightarrow T_{m}-0}\left\Vert u\left( t\right) \right\Vert _{L^{2}}=\infty ,
\end{equation*}%
which completes the proof of Theorem \ref{Thm22}.

At the end of this section, we mention the strategy of the proof of Proposition \ref{Prop32}.
We apply a test-function method used by Zhang \cite{Zha99, Zha01} to NLS (\ref{eq11}).
By using some test-functions and space-time sets cleverly,
he obtained some blow-up results for nonlinear parabolic equations (see \cite{Zha99}).
By the same method, he also proved a blow-up result for the nonlinear damped wave equation:
\begin{equation*}
	\left\{\begin{array}{l}
		v_{tt}-\Delta v+v_t=\vert v\vert ^{p},\quad(t,x) \in \mathbf{R}\times\mathbf{R}^n, \\ 
		v(0,x)=v_0(x), v_t(0,x)=v_1(x),\quad x\in\mathbf{R}^n,
	\end{array}\right.
\end{equation*}
where 
$1<p\leq 1+2/n,$ $v=v\left( t,x\right)$
is a real-valued unknown function,
$v_{0}\left( x\right) $ and $v_{1}\left( x\right) $
are compactly supported given functions (see \cite{Zha01}). However, since this
method needs a positivity of the nonlinear term $\left\vert v\right\vert ^{p}$,
it can not be applicable to NLS (\ref{eq11}) directly,
because solutions for NLS are generally complex-valued and the
constant $\lambda $ in front of the nonlinearity is a complex number.
To overcome these difficulties,
we make a little modification to this method by introducing an appropriate positive function
(see (\ref{3.2})) related to
$\lambda \left\vert u\right\vert ^{p}$.

For the nonlinear heat equation and the damped wave equation with the
same type nonlinearity as $\left\vert u\right\vert ^{p}$,
it is well known that the exponent $p=1+2/n$,
which is often referred to as the \textquotedblleft Fujita exponent\textquotedblright ,
is the threshold between the small data global existence and blow-up of solutions
(see \cite{Lev90, DeLe00, ToYo01} and the references therein).

\section{Proof of Proposition \protect\ref{Prop32}}

In this section, we give a proof of Proposition \ref{Prop32}.
For simplicity, we write $A\simeq B$ if there exist some positive constants
$C_{1},C_{2}>0$ such that $C_{1}B\leq A\leq C_{2}B$ and we also use
$A\lesssim B$ if there exists a positive constant $C>0$ such that $A\leq CB.$

\begin{proof}
First we introduce two cut-off functions
$\eta =\eta (t) \in C_0^{\infty }([0,\infty))$
and
$\phi=\phi(x)\in C_0^{\infty }(\mathbf{R}^n)$
such that
$0\leq\eta ,\phi\leq 1$,
\begin{equation*}
	\eta(t)\equiv\left\{\begin{array}{ll}
		1&\text{ if }t\leq 1/2\\ 
		0&\text{ if }t\geq 1
	\end{array}\right. ,\quad
	\phi(x)\equiv\left\{\begin{array}{ll}
		1&\text{ if }\vert x\vert \leq 1/2\\ 
		0&\text{ if }\vert x\vert \geq 1%
	\end{array}\right. .
\end{equation*}%
Furthermore, it is possible to take $\phi $ satisfying the inequality
\begin{equation}
\label{eq41}
	\frac{\vert( \nabla \phi ) ( x) \vert ^{2}}{\phi (x) }
	\leq C\text{ \ for }\left\vert x\right\vert\leq 1, 
\end{equation}%
with some constant $C$ independent of $x$.
Let $R>0$ be large parameter.
Using the above cut-off functions, we also put three cut-off functions
dependent on $R:$%
\begin{gather}
	\eta _{R}( t) \equiv \eta \left( \frac{t}{R^{2}}\right) \text{ \ for }t\in \mathbf{R},\text{ }\phi _{R}
	\left( x\right) \equiv \phi \left( 
\frac{x}{R}\right) \text{ \ for }x\in \mathbf{R}^{n}, \notag \\
\psi _{R}\left( t,x\right) \equiv \eta _{R}\left( t\right) \phi _{R}\left(
x\right) \text{ \ for }\left( t,x\right) \in \left[ 0,\infty \right) \times 
\mathbf{R}^{n}.  \label{eq42}
\end{gather}%
Let
$B_{R}\equiv \left\{ x\in \mathbf{R}^{n};\text{ }\left\vert x\right\vert \leq R\right\} $
be a ball at the origin. We
also define the time-space set $Q_{R}\equiv [0,R^2]\times B_{R}.$
We note that $Q_R$ includes the support of $\psi_R$.
Denote $q\equiv p/(p-1)\in \left[ 1+n/2,\infty \right) .$
We consider the case $\lambda _{1}>0$ and $\lambda_1\int f_2dx<0$ only,
since the other cases can be treated almost in the same way (see Remark \ref{Remark 3}).
In this case, we may assume $f_{2}\in L^{1}$ and $\int_{\mathbf{R}^{n}}f_{2}\left( x\right) dx<0$
by the assumption (\ref{eq22}).
We define a positive function of $R$ by
\begin{equation}
	I_{R}\equiv\, {\rm Re}\, \int_{Q_{R}}\lambda\vert u\vert^{p}\psi_R^qdxdt.\label{3.2}
\end{equation}%
We note that $\psi _{R}^{q}\in C_{0}^{2}([ 0,R^2+1)\times\mathbf{R}^{n})$.
Since $u$ is a global weak solution of (\ref{eq11}) (see Definition \ref{Def31}),
we can use the identity (\ref{eq31}) with $T=R^{2}+1$ and have
\begin{eqnarray}
I_{R} &=&\int_{B_{R}}f_{2}\left( x\right) \phi _{R}^{q}\left(
x\right) dx+q\int_{Q_{R}}\left( \, {\rm Im}\, u\right) \psi
_{R}^{q-1}\partial _{t}\left( \psi _{R}\right) dxdt \notag \\
&&+\int_{Q_{R}}\left( \, {\rm Re}\, u\right) \Delta \left( \psi
_{R}^{q}\right) dxdt.   \label{3.3}
\end{eqnarray}%
By the assumption on $f_2$, the first term of the right hand side of (\ref%
{3.3}) is negative for sufficiently large $R>0.$
In fact, by $f_2\in L^1$ and Lebesgue's convergence theorem,
there exists $R_{1}>0$ such that for any $R>R_{1},$%
\begin{equation*}
\int_{B_{R}}f_{2}\left( x\right) \phi _{R}^{q}\left( x\right) dx<0.
\label{3.3a}
\end{equation*}
Thus, we have for $R>R_{1}$, 
\begin{eqnarray}
I_{R} &<&q\int_{Q_{R}}\left( \, {\rm Im}\, u\right) \psi
_{R}^{q-1}\left( \partial _{t}\psi _{R}\right) dxdt+\int_{Q%
_{R}}\left( \, {\rm Re}\, u\right) \Delta \left( \psi _{R}^{q}\right) dxdt  \notag
\\
&\lesssim &\int_{Q_{R}}\left\vert u\right\vert \psi
_{R}^{q-1}\left\vert \partial _{t}\left( \psi _{R}\right) \right\vert
dxdt+\int_{Q_{R}}\left\vert u\right\vert \left\vert \Delta \left(
\psi _{R}^{q}\right) \right\vert dxdt  \notag \\
&\equiv &J_{1,R}+J_{2,R}.  \label{3.5}
\end{eqnarray}%
First we will estimate $J_{1,R}.$ By a simple calculation, we get%
\begin{equation*}
\partial _{t}\psi _{R}\left( t,x\right) =\frac{1}{R^{2}}\phi _{R}\left(
x\right) \left( \partial _{t}\eta \right) \left( \frac{t}{R^{2}}\right) .
\end{equation*}%
By noting $\partial _{t}\eta \left( t\right) =0$ if $t\in [0,1/2]$ and the H\"older inequality, we obtain
\begin{eqnarray}
	J_{1,R}&\lesssim &\frac{1}{R^2}\int_{R^2/2}^{R^2}\int_{B_R}\vert u\vert\psi_R^{q-1}dxdt\notag\\
	&\lesssim &\frac{1}{R^2}
		\left( \int_{R^{2}/2}^{R^2}\int_{\mathbf{B}_R}\vert u\vert ^p\psi_R^qdxdt\right) ^{1/p}
		\left( \int_{R^2/2}^{R^2}\int_{B_R}dxdt\right) ^{1/q} \notag\\
	&\simeq &I_{1,R}^{1/p}R^{(n+2-2q)/q},\label{3.6}
\end{eqnarray}%
where%
\begin{equation*}
	I_{1,R}\equiv \, {\rm Re}\, \int_{R^2/2}^{R^2}\int_{B_R}\lambda \vert u\vert ^p\psi_R^qdxdt.
\end{equation*}%
We note that $n+2-2q\leq 0$, since $1<p\leq 1+2/n$.
Next we consider $J_{2,R}.$
By a direct computation, we have 
\begin{eqnarray*}
\Delta \left( \psi _{R}^{q}\right)  &=&\frac{1}{R^{2}}q\left( q-1\right)
\eta _{R}^{q}\left( t\right) \phi _{R}^{q-2}\left( x\right) \left\vert
\nabla \phi \right\vert ^{2}\left( \frac{x}{R}\right)   \notag \\
&&+\frac{1}{R^{2}}q\eta _{R}^{q}\left( t\right) \phi _{R}^{q-1}\left(
x\right) \left( \Delta \phi \right) \left( \frac{x}{R}\right) .  \label{3.6a}
\end{eqnarray*}%
Using this and (\ref{eq41}), in the same manner as above, we obtain
\begin{align}
	J_{2,R}&\lesssim\frac{1}{R^2}\left( \int_0^{R^2}\int_{B_R\setminus B_{R/2}}
		\vert u\vert ^p\psi_R^qdxdt\right)^{1/p}
		\left( \int_0^{R^2}\int_{B_R\setminus B_{R/2}}dxdt\right)^{1/q}\notag\\
	&\simeq I_{2,R}^{1/p}R^{(n+2-2q)/q},\label{3.7}
\end{align}%
where we put 
\begin{equation*}
I_{2,R}\equiv \, {\rm Re}\, \int_{0}^{R^{2}}\int_{B_{R}\backslash 
B_{R/2}}\lambda \left\vert u\right\vert ^{p}\psi
_{R}^{q}dxdt.
\end{equation*}%
By combining (\ref{3.5}), (\ref{3.6}) and (\ref{3.7}), we have%
\begin{equation}
I_{R}\lesssim \left( I_{1,R}^{1/p}+I_{2,R}^{1/p}\right) R^{(n+2-2q)/q},  \label{3.8}
\end{equation}%
for $R>R_{1}.$
Since it is clear that 
$I_{j,R}\leq I_{R} \,( j=1,2)$,
we obtain
\begin{equation}
I_{R}\lesssim R^{n+2-2q}\leq C,  \label{3.9}
\end{equation}%
with some constant $C$ independent of $R$, since $n+2-2q\leq 0$.
Here we note that only in the critical case $p=1+2/n,$ the identity $n+2-2q=0$ holds.
By (\ref{3.9}) and letting $R\rightarrow +\infty$, we have
\begin{equation*}
\, {\rm Re}\, \int_{\left[ 0,\infty \right) \times \mathbf{R}^{n}}\lambda
\left\vert u\right\vert ^{p}dtdx<\infty ,  \label{3.10}
\end{equation*}%
that is, $u\in L^{p}\left(\left[ 0,\infty \right) \times \mathbf{R}^{n}\right) .$
Noting this and the integral region of $I_{1,R}$ and $I_{2,R}$, we have
\begin{equation}
\lim_{R\rightarrow +\infty }I_{j,R}=0,\text{ \ for }j=1,2.  \label{3.10c}
\end{equation}%
Therefore by the inequality (\ref{3.8}) and (\ref{3.10c}), we get 
\begin{equation*}
\lim_{R\rightarrow +\infty }I_{R}=0,  \label{3.11}
\end{equation*}
which implies $u=0$.
This completes the proof.
\end{proof}

\begin{remark}
\label{Remark 3}
In the different cases from $\lambda_1>0$,
putting
$$
	I_R\equiv \left\{\begin{array}{cl}
		\displaystyle-\, {\rm Re}\, \int_{Q_R}\lambda\vert u\vert ^p\psi_R^qdxdt
			&\text{if}\ \lambda_1<0, \ \lambda_1\int f_2dx<0,\\[8pt]
		\displaystyle\, {\rm Im}\, \int_{Q_R}\lambda\vert u\vert ^p\psi_R^qdxdt
			&\text{if}\ \lambda_2>0, \ \lambda_2\int f_1dx>0,\\[8pt]
		\displaystyle-\, {\rm Im}\, \int_{Q_R}\lambda\vert u\vert ^p\psi_R^qdxdt
			&\text{if}\ \lambda_2<0, \ \lambda_2\int f_1dx>0,
		\end{array}\right.
$$
we can prove the same conclusion in the same manner as above.
\end{remark}

\section{\label{S4} Appendix}
 
In this section, we give a proof of Proposition \ref{Prop31}.
The main difficulty of the proof lies in the fact
that if $p$ is close to $1,$ then the nonlinear term $\left\vert
u\right\vert ^{p}$ does not have twice differentiability with
respect to space variables. To avoid differentiating twice, we use
appropriate changing variables and differentiate with regard to time
variable (see (\ref{1.12b})). As the result, we can derive an $H%
^{2} $-estimate (see also \cite{Caze03}).

We first recall the well-known Strichartz
estimates for the Schr\"{o}dinger equation (see \cite{Ya87}).

Let%
\begin{equation*}
	\left\{ \begin{array}{c}
		2\leq \rho_j <2n/(n-2)\text{ \ if }n\geq 3 \\ 
		2\leq \rho_j <\infty \text{ \ if }n=2 \\ 
		2\leq \rho_j \leq \infty \text{ \ if }n=1%
	\end{array}\right. \text{ and }\ \frac{2}{r_j}=\frac{n}{2}-\frac{n}{\rho_j} \ (j=1,2) .
\end{equation*}%
Then the following estimates hold:

\begin{lemma}
\label{Lemma 1} For any time interval $I$, the estimates 
\begin{align}
	\Vert U(t) f\Vert _{L_{t}^{r_1}( I;L_x^{\rho_1 }) }
	&\lesssim\Vert f\Vert _{L^2},\notag\\
	\left\Vert \int_{0}^{t}U( t-s) Fds\right\Vert _{L_{t}^{r_1}( I;L_{x}^{\rho_1 }) }
	&\lesssim \Vert F\Vert _{L_{t}^{r_2^{\prime }}( I;L_x^{\rho_2 ^{^{\prime }}})}
	\label{eq51}
\end{align}
are true, where $r_2^{\prime }=r_2/(r_2-1)$ and $\rho_2^{\prime}=\rho_2/(\rho_2-1)$.
\end{lemma}

Now we give a proof of Proposition \ref{Prop31}.
Denote the nonlinear term by
$F\left( u\right) =\lambda \left\vert u\right\vert ^{p}$
and the time interval by $I=\left[ 0,T\right) $ for simplicity.

\begin{proof}
Let $T>0, \rho=p+1, 2/r=n/2-n/\rho$,
$\psi\in C^2_0([0,T)\times\mathbf{R}^n)$
and let
$u$ be an $L^2$-solution of (\ref{eq21}) on $[0,T)$.
It is easy to see that
$u\in L_{loc}^{p}(\left[ 0,T) \times \mathbf{R}^{n}\right)$.
We decompose $u$ into
$u=u_{1}+u_{2}$,
where $u_{1}\equiv U\left( t\right) f$ is the homogeneous part and
\begin{equation*}
	u_{2}\equiv -i\int_{0}^{t}U\left( t-s\right) F\left( u\right) ds
\end{equation*}
is the inhomogeneous one. The homogeneous part $u_{1}$ can be treated
easily. In fact, by a standard density argument, we can obtain the identity
\begin{equation*}
	\int_{I\times\mathbf{R}^n}u_1( -i\partial_t\psi+\Delta\psi)dxdt
	=i\int_{\mathbf{R}^n}f(x)\psi(0,x)dx.\label{eq52}
\end{equation*}
Thus, it suffices to prove
\begin{equation}
	\int_{I\times\mathbf{R}^n}u_2( -i\partial _{t}\psi+\Delta\psi)dxdt
	=\int_{I\times\mathbf{R}^n}F(u)\psi dxdt,\label{eq53}
\end{equation}
which must be dealt with somewhat carefully because of involving the non-smooth nonlinearity
$\vert u\vert ^{p}$.
we split the left-hand-side of (\ref{eq53}) as
\begin{equation}
\label{eq54}
	-i\int_{I\times\mathbf{R}^n}u_2(\partial_t\psi)dxdt
	+\int_{I\times\mathbf{R}^n}u_2\Delta\psi dxdt\\
	\equiv K_1+K_2.
\end{equation}
Hereafter we use the notation
$L_{t}^{r}L_{x}^{\rho }\equiv L_{t}^{r}\left( I;L_{x}^{\rho }\right) $
for simplicity.
Since $u\in L_{t}^{r}L_{x}^{\rho }$
and $C_{0}^{\infty }\left( I\times \mathbf{R}^{n}\right) $
is dense in $L_{t}^{r}L_{x}^{\rho }$,
there exists a sequence
$\{ u_k\}_{k\in\mathbf{N}}\subset C_{0}^{\infty }\left( I\times \mathbf{R}^{n}\right) $
such that
\begin{equation}
\label{eq55}
	\lim_{k\to\infty}\Vert u_k-u\Vert _{L_t^r( I;L_x^{\rho })}=0.
\end{equation}
We also introduce an approximate function sequence
$\left\{ u_{2,k}\right\} _{k\in\mathbf{N}}$
to the inhomogeneous part $u_{2}$,
whose component is given by
\begin{equation*}
	u_{2,k}\equiv -i\int_{0}^{t}U(t-s) F(u_k) ds.
\end{equation*}
Let
$\alpha \equiv \frac{n}{4}\left( 1+\frac{4}{n}-p\right) >0.$
By the Strichartz estimate (\ref{eq51}) and the H\"{o}lder inequality with
$\frac{1}{\rho ^{^{\prime }}}=\frac{1}{\rho }+\frac{p-1}{p+1}$
and
$\frac{1}{r^{^{\prime }}}=\frac{1}{r}+\frac{p-1}{r}+\alpha ,$
we can estimate
\begin{eqnarray}
\left\Vert u_{2}-u_{2,k}\right\Vert _{L_{t}^{\infty }L%
_{x}^{2}} &\lesssim &\left\Vert \left\vert u\right\vert ^{p}-\left\vert
u_{k}\right\vert ^{p}\right\Vert _{L_{t}^{r^{^{\prime }}}L%
_{x}^{\rho ^{^{\prime }}}}  \notag \\
&\lesssim &\left\Vert \left( \left\Vert u\right\Vert _{L_{x}^{\rho
}}^{p-1}+\left\Vert u_{k}\right\Vert _{L_{x}^{\rho }}^{p-1}\right)
\left\Vert u-u_{k}\right\Vert _{L_{x}^{\rho }}\right\Vert _{L_{t}^{r^{^{\prime }}}}  \notag \\
&\lesssim &T^{\alpha }\left( \left\Vert u\right\Vert _{L_{t}^{r}%
L_{x}^{\rho }}^{p-1}+\left\Vert u_{k}\right\Vert _{L%
_{t}^{r}L_{x}^{\rho }}^{p-1}\right) \left\Vert u-u_{k}\right\Vert _{%
L_{t}^{r}L_{x}^{\rho }}.  \label{eq56}
\end{eqnarray}%
By (\ref{eq55}), (\ref{eq56}) , noting $u_{2,k}(0,x)=0$ and integration by parts, we have
\begin{align}
	K_1&=-i\lim_{n\to\infty}\int_{I\times\mathbf{R}^{n}}u_{2,k}(\partial_{t}\psi) dxdt\notag\\
	&=\lim_{n\to\infty }i\int_{I\times \mathbf{R}^{n}}(\partial _{t}u_{2,k})\psi dxdt.\label{eq57}
\end{align}%
By the almost same argument as in (\ref{eq56}),
we find that $u_{2,k}\in C\left( I;H^{1}\right) $
and there exists a time derivative
$\partial _{t}u_{2,k}\in C\left( I;H^{-1}\right) $
such that the identity
\begin{equation}
\label{eq58}
	\partial _{t}u_{2,k}=i\Delta u_{2,k}-iF\left( u_{k}\right)
\end{equation}
is valid.
From this identity, we can show that $\Delta u_{2,k}\in C\left( I;L^{2}\right)$.
In fact, changing variables with $t-s=s^{^{\prime }},$ we have
\begin{align}
	\partial_tu_{2,n}(t)
	&=-i\partial_t\int_{0}^{t}U(s^{\prime})F(u_k)(t-s^{\prime})ds^{\prime}\notag\\
	&=-iU(t)F(u_k)(0)-i\int_0^tU(s)\partial_t(F(u_k))(t-s)ds.\label{1.12b}
\end{align}%
Applying the Strichartz estimate (\ref{eq51}) to (\ref{1.12b}), we have
\begin{equation}
\left\Vert \partial _{t}u_{2,k}\left( t\right) \right\Vert _{L%
^{2}}\lesssim \left\Vert u_{k}\left( 0\right) \right\Vert _{L%
^{2p}}^{p}+\left\Vert \partial _{t}\left( F\left( u_{k}\right) \right)
\right\Vert _{L_{t}^{r^{^{\prime }}}L_{x}^{\rho ^{^{\prime}}}}.  \label{1.7a}
\end{equation}%
By the same way as in (\ref{eq56}), we also have
\begin{equation}
\left\Vert \partial _{t}\left( F\left( u_{k}\right) \right) \right\Vert _{%
L_{t}^{r^{^{\prime }}}L_{x}^{\rho ^{^{\prime }}}}\lesssim
T^{\alpha }\left\Vert u_{k}\right\Vert _{L_{t}^{r}L%
_{x}^{\rho }}^{p-1}\left\Vert \partial _{t}u_{k}\right\Vert _{L%
_{t}^{r}L_{x}^{\rho }}\label{1.7b}
\end{equation}%
and the right-hand-side is finite due to
$u_k\in C_0^{\infty}( I\times\mathbf{R}^{n})$.
Therefore by combining (\ref{1.7a})-(\ref{1.7b}), we obtain%
\begin{equation*}
\left\Vert \partial _{t}u_{2,k}\left( t\right) \right\Vert _{L%
^{2}}\lesssim \left\Vert u_{k}\right\Vert _{L_{t}^{\infty }L_{x}^{2p}}^{p}
+T^{\alpha }\left\Vert u_{k}\right\Vert _{L_{t}^{r}%
L_{x}^{\rho }}^{p-1}\left\Vert \partial _{t}u_{k}\right\Vert _{%
L_{t}^{r}L_{x}^{\rho }}<\infty  \label{1.7c}
\end{equation*}%
for any $k\in \mathbf{N}$,
from which we can see $\partial _{t}u_{2,k}\in C\left( I;L^{2}\right)$.
Thus by the equation (\ref{eq58}) again, we also find
$u_{2,k}\in C\left( I;H^{2}\right)$
for any $k\in \mathbf{N}$.
Therefore we have the identity 
\begin{equation}
\left( \Delta u_{2,k},\psi \right) _{L_{x}^{2}}=\left(
u_{2,k},\Delta \psi \right) _{L_{x}^{2}}.  \label{1.8b}
\end{equation}%
Thus by combining the identities (\ref{eq57}), (\ref{eq58}) and (\ref{1.8b}), we obtain
\begin{align}
	K_1&=\lim_{k\to\infty }\left( \int_{I\times\mathbf{R}^n}F( u_k)\psi dxdt
		-\int_{I\times \mathbf{R}^n}u_{2,k}\Delta \psi dxdt\right)\notag\\
	&=\int_{I\times \mathbf{R}^{n}}F(u)\psi dxdt-K_2.\label{1.8c}
\end{align}
In fact, by the same way as in (\ref{eq56}), we obtain
\begin{equation*}
\left\vert \int_{I\times \mathbf{R}^{n}}\left( F\left(
u_{k}\right) -F\left( u\right) \right) \psi dxdt\right\vert  \label{1.13a}
\lesssim T^{\alpha }\left( \left\Vert u_{k}\right\Vert _{L_{t}^{r}%
L_{x}^{\rho }}^{p-1}+\left\Vert u\right\Vert _{L_{t}^{r}%
L_{x}^{\rho }}^{p-1}\right) \left\Vert u_{k}-u\right\Vert _{L_{t}^{r}L_{x}^{\rho }}
\left\Vert \psi \right\Vert _{L%
_{t}^{r}L_{x}^{\rho }}
\end{equation*}%
and
\begin{equation*}
\left\vert \int_{I\times \mathbf{R}^{n}}\left(
u_{2,k}-u_{2}\right) \Delta \psi dxdt\right\vert  \label{1.13}
\lesssim T\Vert u_{2,k}-u_2\Vert_{L_t^{\infty}L_x^2}\Vert \Delta \psi \Vert_{L_{t}^{\infty }L_{x}^{2}}.
\end{equation*}%
Therefore, combining (\ref{eq54}) and (\ref{1.8c}), we obtain (\ref{eq53}). This completes the proof.
\end{proof}

\emph{Acknowledgments}. The authors would like to express deep gratitude to
Professor Nakao Hayashi and Professor Tatsuo Nishitani for their helpful
advice, comments and constant encouragements. The authors would like also to
thank Professor Soichiro Katayama and Professor Hideaki Sunagawa for many
useful sugesstions and comments.

\end{document}